\documentclass[12pt]{amsart}
\usepackage{amscd}
\def\NZQ{\Bbb}               % the font for N,Z,Q,R,C
\def\NN{{\NZQ N}}

\def\RR{{\NZQ R}}

\def\D{{\Delta}}
\def\frk{\frak}               % font for "Fraktur"

\def\Phi{{\frk n}}
\def\Phi{{\frk N}}

\def\a{{\bold a}}
\def\b{{\bold b}}
\def\c{{\bold c}}
\def\d{{\bold d}}

\def\1{{\mathbf 1}}
\def\0{{\mathbf 0}}
\def\l{{\lambda}}
\def\opn#1#2{\def#1{\operatorname{#2}}} % to make operators
\opn\relint{relint} 
\opn\height{ht}

\newtheorem{Theorem}{Theorem}[section]
\newtheorem{Lemma}[Theorem]{Lemma}
\newtheorem{Corollary}[Theorem]{Corollary}

\newtheorem{Remark}[Theorem]{Remark}

\newtheorem{Problem}[Theorem]{Problem}

\begin{document}

\title{Integral closures of monomial ideals\\  and Fulkersonian hypergraphs}
\author{Ngo Viet Trung}
\address{Institute of Mathematics, 18 Hoang Quoc Viet, 10307 Hanoi, Vietnam}
\email{nvtrung@math.ac.vn}
\begin{abstract}
We prove that the integral closures of the powers of a squarefree monomial ideal $I$ equal the symbolic powers  if and only if $I$ is the edge ideal of a Fulkersonian hypergraph. \vskip 0.7cm
\centerline{\it Dedicated to Do Long Van on the occasion of his sixtyfifth birthday}
\end{abstract} 
\maketitle

\section*{Introduction}

Let $V$ be a finite set. A hypergraph $\D$ on $V$  is a family of subsets of $V$. The elements of $V$ and $\D$ are called the vertices and the edges of $\D$, respectively. We call $\D$ a simple hypergraph if there are no inclusions between the edges of $\D$.  

Assume that $V = \{1,...,n\}$ and let $R = K[x_1,...,x_n]$ be a polynomial ring over a field $K$. The edge ideal $I(\D)$ of $\D$ in $R$ is the ideal generated by all monomials of the form $\prod_{i \in F}x_i$ with $F \in \D$. By this way we obtain an one-to-one correspondence between simple hypergraphs and squarefree monomials.

It is showed  \cite{HHTZ} (and implicitly in \cite{GVV}) that the symbolic powers of $I(\D)$ coincide with the ordinary powers of $I(\D)$ if and only if $\D$ is a Mengerian hypergraph, which is defined by a min-max equation in Integer Linear Programming. A natural generalization of the Mengerian hypergraph is the Fulkersonian hypergraph which is defined by the integrality of the blocking polyhedron. Mengerian and Fulkersonian hypergraphs belong to a variety of hypergraphs which generalize bipartite graphs and trees in Graph Theory \cite{Be} \cite{Du}.  They frequently arise in the polyhedral approach of combinatorial optimization problems.

The aim of this note is to show that the symbolic powers of $I(\D)$ coincide with the integral closure of the ordinary powers of $I(\D)$ if and only if $\D$ is a Fulkersonian hypergraph. We will follow the approach of \cite{HHT} \cite{HHTZ} which describes the symbolic powers of squarefree monomials by means of the vertex covers of hypergraphs. This approach will be presented in Section 1. The above characterization of the integral closure of the ordinary powers of squarefree monomials ideals will be proved in Section 2.

\section{Vertex covers and symbolic powers}

Let $\D$ be a simple hypergraph on $V =  \{1,...,n\}$.
For every edge $F \in \D$ we denote by $P_F$ the ideal $(x_i|\ i \in F)$ in the polynomial ring $R = K[x_1,...,x_n]$. Let
$$I^*(\D) := \bigcap_{F \in \D} P_F.$$
Then $I^*(\D)$ is a squarefree monomial ideal in $R$.   It is clear that every squarefree monomial ideal can be viewed as an ideal of the form $I^*(\D)$.

A subset $C$ of $V$ is called a vertex cover of $\D$ if it meets every edge. 
Let $\D^*$ denote the hypergraph of the minimal vertex covers of $\D$. 
This hypergraph is known under the name transversal \cite{Be} or blocker  \cite{Du}.
It is well-known that 
$I^*(\D) = I(\D^*).$
For this reason we call $I^*(\D)$ the vertex cover ideal of $\D$.

Viewing a squarefree monomial ideal $I$ as the vertex cover ideal of a hypergraph is  suited for the study of the symbolic powers of $I$.  
If $I = I^*(\D)$, then the $k$-th symbolic power of $I$ is the ideal
$$I^{(k)} =   \bigcap_{F \in \D} P_F^k.$$
The monomials of $I^{(k)}$ can be described by means of $\D$ as follows \cite{HHT}.

Let $\c = (c_1,...,c_n)$ be an arbitrary integral vector in $\NN^n$. We may think of $\c$ as a multiset consisting of $c_i$ copies of $i$ for $i = 1,...,n$. Thus, a subset $C \subseteq V$ corresponds to an (0,1)-vector $\c$ with $c_i = 1$ if $i \in C$ and $c_i = 0$ if $i \not\in C$, and $C$ is a vertex cover of $\D$ if $\sum_{i \in F}c_i \ge 1$ for all $F \in \D$.  
For this reason, we call  $\c$  a vertex cover of order $k$ of $\D$  if $\sum_{i \in F}c_i \ge k$ for all $F \in \D$. 
Let $x^\c$ denote the monomial $x_1^{c_1} \cdots x_n^{c_n}$. It is obvious that $x^\c \in P_F$  iff $\sum_{i \in F}c_i \ge k$. Therefore, $x^\c \in I^{(k)}$ iff $\c$ is a vertex cover of order $k$. In particular, $x^\c \in I$ iff $\c$ is a vertex cover of order $1$.

Let $F_1,...,F_m$ be the edges of $\D$. We may think of  $\D$ as an $n \times m$ matrix $M = (e_{ij})$ with $e_{ij} =  1$ if $i \in F_j$ and $e_{ij} = 0$ if $i \not\in F_j$. One calls $M$ the incidence matrix of $\D$. 
Since the columns of $M$ are the integral vectors of $F_1,...,F_m$, an integral vector $\c \in \NN^n$ is a vertex cover of order $k$ of $\D$ iff $M^T \cdot \c \ge k\1$, where $\1 $ denote the vector $(1,...,1)$ of  $\NN^m$.

By the above characterization of monomials of symbolic powers we have  $I^{(k)} = I^k$ if every vertex cover $\c$ of order $k$ can be decomposed as as sum of $k$ vertex cover of order 1 of $\D$.

Every integral vector $\c \in \NN^r$ is a vertex cover of some order $k \ge 0$. The minimum order of $\c$ is the number $o(\c) := \min\{\sum_{i \in F}c_i|\ F \in \D\}.$ 
Let $\sigma(\c)$ denote the maximum number $k$ such that $\c$ can be decomposed as a sum of $k$ vertex cover of order 1. Then $I^{(k)} = I^k$ for all $k \ge 1$ if and only if $o(\c) = \tau(\c)$ for 
every integral vector $\c \in \NN^r$.

Using the incidence matrix of the hypergraph of minimal vertex covers one can characterize the numbers $o(\c)$ and $\tau(\c)$ as follows.

\begin{Lemma} \cite[Lemma 1.3]{HHTZ}
Let $M$ be the incidence matrix of the  hypergraph $\D^*$ of the minimal vertex covers of $\D$. Then
\par {\rm (i) }  $o(\c) = \min\{\a \cdot \c|\  \a \in \NN^n,\  M^T\cdot\a \ge \1 \}$,
\par {\rm (ii)}  $\sigma(\c) = \max\{\b \cdot \1 |\ \b \in \NN^m,\ M\cdot\b \le \c\}$.
\end{Lemma}

Let $M$ now be the incidence matrix of a hypergraph $\D$. One calls $\D$ a Mengerian hypergraph \cite{Be} \cite {Du} (or having the max-flow min-cut property \cite{GVV}) if
$$\min\{\a \cdot \c|\  \a \in \NN^n,\  M^T\cdot\a \ge \1 \}
= \max\{\b \cdot \1 |\ \b \in \NN^m,\ M \cdot\b \le \c\}.$$

Since $I(\D) = I^*(\D^*)$,  switching the role of $\D$ and $\D^*$ in the above observations  we immediately obtain the following criterion for the equality of ordinary and symbolic powers of a squarefree monomial ideal.

\begin{Theorem}  \label{Menger} \cite[Corollary 1.6]{HHTZ}  Let $I = I(\D)$. Then $I^{(k)} = I^k$ for all $k \ge 1$ if and only if $\D$ is a Mengerian hypergraph.
\end{Theorem}

\begin{Remark}
{\rm In general, $\D^*$ needn't to be a Mengerian hypergraph if $\D$ is a Mengerian hypergraph (see e.g. \cite[Example 2.8]{HHTZ}).}
\end{Remark}
  
It should be noticed that $\min\{\a \cdot \1|\  \a \in \NN^n,\  M^T\cdot\a \ge \1 \}$ is the 
minimum number of vertices of vertex covers  and $\max\{\b \cdot \1 |\ \b \in \NN^m,\ M \cdot\b \le \1\}$ is the maximum number of disjoint edges of $\D$.
If these numbers are equal, one says that $\D$ has the K\"onig property \cite{Be}, \cite{Du}. 
This is a typical property of trees and bipartite graphs.

\section{Fulkersonian hypergraphs}

Let $\D$ be a simple graph of  $m$ edges on $n$ vertices. Let $M$ be the incidence matrix of $\D$. 
By the duality in Linear Programming we   have
$$ \min\{\a \cdot \c|\  \a \in \RR_+^n,\  M^T\cdot\a \ge \1 \} 
= \max\{\b \cdot \1 |\ \b \in \RR_+^m,\ M \cdot\b \le \c\},$$
where $\RR_+$ denote the set of non-negative reel numbers. This implies
$$\min\{\a \cdot \c|\  \a \in \NN^n,\  M^T\cdot\a \ge \1 \}
\le   \max\{\b \cdot \1 |\ \b \in \NN^m,\ M \cdot\b \le \c\}.$$
If equality holds above,  we obtain
\begin{align*}
\min\{\a \cdot \c|\  \a \in \RR_+^n,\  M^T\cdot\a \ge \1 \} & =   \min\{\a \cdot \c|\  \a \in \NN^n,\  M^T\cdot\a \ge \1 \},\\
\max\{\b \cdot \1 |\ \b \in \RR_+^m,\ M \cdot\b \le \c\} & = \max\{\b \cdot \1 |\ \b \in \NN^m,\ M \cdot\b \le \c\}.
\end{align*}
In this case, the two optimization problems  on the left sides have integral optimal solutions. 

For the min problem, this condition is closely related to the integrality of the polyhedron:
$$Q(\D)  := \{\a \in \RR_+^n|\  M^T\cdot\a \ge \1\}.$$ 
This polyhedron is usually called the blocking polyhedrone of $\D$ \cite{Du}.
Notice that an integral vector $\c  \in \NN^n$ is a vertex cover of order $1$ of $\D$ iff  $\c \in Q(\D)$. 

\begin{Lemma} \label{Hoffman} {\rm (see e.g. \cite[Lemma 1, p.~203]{Be})}
$\min\{\a \cdot \c|\  \a \in \RR^n,\  M^T\cdot\a \ge \1 \}$ is an integer for all $\c \in \NN^n$ if and only if $Q(\D)$ only has integral extremal points.
\end{Lemma}

One calls  $\D$ a Fulkersonian hypergraph \cite{Du} (or paranormal \cite{Be}) if  $Q(\D)$  only has integral extremal points.  By the above observation and Lemma \ref{Hoffman}, Fulkersonian hypergraphs are generalizations of Mengerian hypergraphs. 

Unlike the Mengerian property, the Fulkersonian property is preserved by passing to the hypergraph of minimal vertex covers.

\begin{Lemma} \label{Lehman} {\rm (see e.g. \cite[Corollary, p.~210]{Be})}
$\D$ is Fulkersonian if and only if $\D^*$ is Fulkersonian.
\end{Lemma}

We shall see that Fulkersonian hypergraphs can be used to study the integral closures of powers of monomial ideals.

Let $I$ be an arbitrary monomial ideals. Let $\overline{I}$ denote the integral closure of $I$. It is easy to see  that $\overline{I}$  is the monomial ideal generated by all monomial $f$ such that $f^p \in I^p$ for some $p \ge 1$. We say that $I$ is an integrally closed ideal if $\overline{I} = I$.

It is well known that powers of ideals generated by variables are integrally closed. Since the intersection of integrally closed ideals is again an integrally closed ideal, symbolic powers of squarefree monomial ideals are integrally closed. From this it follows that $\overline{I^k} \subseteq I^{(k)}$ for all $k \ge 0$ if $I$ is a squarefree monomial ideal.

\begin{Theorem} \label{Fulkerson-1}
Let  $I = I^*(\D)$. Then $\overline{I^k} = I^{(k)}$ for all $k \ge 1$ if and only if $\D$ is a Fulkersonian hypergraph.
\end{Theorem}

\begin{proof}
Assume that $Q(\D)$ is integral with integral vertices $\a_1,...,\a_r$. We have to show that every monomial $x^\c \in I^{(k)}$ belongs to $I^k$.
As we have seen in Section 1, $\c$ is a vertex cover of order $k$ of $\D$. This means $M^T\cdot\c \ge k\1$. Therefore 
$\frac{1}{k}\c  \in Q(\D)$. Hence there are rational numbers $\l_1,..,\l_r \ge 0$ with $\l_1 + \cdots + \l_r = 1$ such that
$$\frac{1}{k}\c  = \l_1\a_1 + \cdots + \l_s\a_r + \b$$
for some rational vector $\b \in \RR_+^n$. Let $p$ be the least common multiple of the denominators of $\l_1,...,\l_r$ and the components of $\b$. Then 
$$p\c = kp\l_1\a_1 + \cdots + kp\l_r\a_r + kp\b$$
is a sum of $kp$ integral vectors $\a_1,...,\a_r$  in $Q(\D)$ and the integral vector $kr\b \in \NN^n$. 
Since $x^{\a_1},...,x^{\a_r} \in I$,  
$$(x^\c)^p = (x^{\a_1})^{kp\l_1} \cdots (x^{\a_r})^{kr\l_r}x^{kp\b} \in I^{kp}.$$
Therefore, $x^\c \in \overline{I^k}$ as required.

Conversely, assume that $I^{(k)} = \overline{I^k}$ for all $k \ge 1$.
Let $\a_1,...,\a_r$ now be the integral vectors corresponding to the minimal vertex covers of $\D$.
Let $P(\D)$ denote the set of all vectors of the form $\l_1\a_1 + \cdots + \l_r\a_r + \b$ with $\mu_1,...,\mu_r \in \RR_+$ and $\b \in \RR_+^n$.
It is obvious that $P(\D) \subseteq Q(\D)$. We will prove that $Q(\D) = P(\D)$, which shows that $\a_1,...,\a_r$ are the extremal points of $Q(\D)$.  

It suffices to show that every rational vector $\a \in Q(\D)$ belongs to $P(\D)$. Let $k$ be the least common multiple of the denominators of the components of $\a$.
Then $M^T\cdot (k\a) \ge k\1$. Hence $x^{k\a} \in I^{(k)} = \overline{I^k}$.  Thus,  there exists an integer $p \ge 1$
such that $x^{pk\a} \in I^{pk}$.  Since $I$ is generated by $x^{\a_1},...,x^{\a_r}$, we have
$x^{pk\a} = x^{\nu_1\a_1}\cdots x^{\nu_r\a_r}x^\d $
for some integral vector $\d \in \NN^n$ and integers $\nu_1,....,\nu_r$ with $\nu_1 + \cdots + \nu_r = pk$. It follows that
$$\a =  \frac{\nu_1}{pk}\c_1 + \cdots + \frac{\nu_r}{pk}\c_r + \frac{1}{pk}\d.$$
Therefore, $\a \in P(\D)$, as desired.
\end{proof}

By Lemma \ref{Lehman}, Theorem \ref{Fulkerson-1} can be reformulated as follows.

\begin{Theorem} \label{Fulkerson-2}
Let  $I = I(\D)$. Then $\overline{I^k} = I^{(k)}$ for all $k \ge 1$ if and only if $\D$ is a Fulkersonian hypergraph.
\end{Theorem}

It is obvious that $I^{(k)} = I^k$ for all $k \ge 1$ iff $I^{(k)} = \overline{I^k}$ and $\overline{I^k} = I^k$ for all $k \ge 1$. Let $R[It] = \bigoplus_{k \ge 0}I^kt^k$ be the Rees algebra of $I$.
It is known that $R[It]$ is normal iff $\overline{I^k} = I^k$ for all $k \ge 1$.
Therefore, combining Theorem \ref{Menger} and Theorem \ref{Fulkerson-2} we obtain the following result of 
Gitler, Valencia and Villarreal \cite[Theorem 3.5]{GVV}.

\begin{Corollary}  \label{GVV}
Let $I = I(\D)$. Then $\D$ is a Mengerian hypergraph if and only if $\D$ is a Fulkersonian hypergraph and $R[It]$ is normal. 
\end{Corollary}

In an earlier paper, Escobar,  Villarreal and  Yoshino showed   that $I^{(k)} =  I^k$ for all $k \ge 1$ iff $\D$ is a Fulkersonian hypergraph and $R[It]$ is normal \cite[Proposition 3.4]{Es}. Combining this result with Corollary \cite{GVV} one can recover Theorem \ref{Menger}.

In view of Corollary \ref{GVV} it is of great interest to study the following 

\begin{Problem}
Let $I = I(\D)$. Can one describe the normality of the Rees algebra $R[It]$ in terms of  $\D$?
\end{Problem}

This problem has been solved for the graph case by Hibi and Ohsugi \cite{HO}, Simis, Vasconcelos and %%@
Villarreal \cite{SVV}.

\end{document}